%% file: artin.tex
\documentclass[12pt]{amsart}
\usepackage{amsmath,amsfonts,amssymb,amscd, epsfig}
\usepackage{latexsym}

\usepackage[
hyperindex=true,pagebackref=true,bookmarks=true,
colorlinks=true,linkcolor=blue,citecolor=red]
{hyperref}
\def\~{{\rm --}} 

\title [A note on Artin's constant]
{A note on Artin's constant} 
\author[Ivan Cherednik]{Ivan Cherednik $^\dag$}

\thanks{$^\dag$  \today\ \ \ Partially supported by NSF grant
DMS--0800642}

\address[I. Cherednik]{Department of Mathematics, UNC
Chapel Hill, North Carolina 27599, USA\\
chered@email.unc.edu}

 \def\~{{\bf --}}
\newcommand{\comment}[1]{}

\input macro

\begin{document}
\maketitle

\renewcommand{\natural}{\wr}


\comment{
We suggest a new approach to Artin's constant that
leads to its representation as an infinite sum divided 
by another infinite sum. The same approach
works well for Stephens' constant and higher rank Artin's
constants. The main results are theoretical but there are 
interesting experimental and computational aspects.
}

\renewcommand{\thesubsection}{\arabic{subsection}}

We suggest a representation of 
Artin's constant, which conjecturally describes the
density of prime $p$ such that ``generic" $g\in \Z$ is 
primitive modulo $p$. Namely, $A=\lim_{N\to \infty} R_k(N)$ 
for $R_k(N)=\frac{\sum p^k\phi(p_i-1)}{\sum p^k(p_i-1)}$, where
the summation is over first $N$ prime $p_i$,\, $k\in \Z_+$.
The classical summation formula is as follows: 
$A=\lim_{N\to \infty}\Si(N)$, where 
$\Si(N)=\frac{1}{N}\sum\frac{\phi(p_i-1)}{p_i-1}$.
The changes needed for arbitrary $g$ are addressed in 
Theorem \ref{AG}, a good exercise in basic analytic and
algebraic number theory. The same 
procedure can be applied to other number-theoretic constants 
like $A$ (see, e.g., \cite{Ni}). In Theorem 2, 
we demonstrate how it works for the Stephens constant and 
for Artin's constants of higher ranks (for the density of 
prime $p$ such that a given set of ``generic" integers generates 
$\Z_p^*$\,). 

The following three features of this approach vs. the
summation formulas are worth noticing. 

1) The {\em restricted summation} suggested by P. Moore to make 
the $\Si$\~formula matching the right heuristic density
for arbitrary $g\in \Z$ gives the desired answer in our approach 
only when $g$ is not a pure (odd) power in $\Z$. Otherwise, 
nontrivial rational multiplicative corrections occur; 
they are calculated in  Theorem \ref{AG}, (ii).

2) The {\em $p^k$\~terms}
in the denominator and numerator of $R_k(N)$
do not influence the limit, which can be 
heuristically associated with switching from primitive
roots in $\Z_p$ to those in $(\Z/(p^{k+1}))^*$. The extra 
$p^k$\~factors disappear (cancel) in the corresponding summation 
formula. When $k=0$, our $R$\~formula for $A$ 
(without restricting the summation) follows from \cite{Pi}.

3) The $R$\~formulas oscillate significantly around $A$ 
(and the other constants). The magnitude of {\em oscillations}   
increases as $k$ grows; see Figure \ref{art0-1}. Representing
$R_k(N)=\sum_{i=1}^N w_i\frac{\phi(p_i-1)}{p_i-1}$, the 
weights $w_i$ change from  $O(\frac{\log N} {N^{k+2}})$ 
for small $i$ to $O(\frac{(\log N)^{k+1}} {N})$ for $i\sim N$. 
Thus, large $\,p\,$ receive greater weights in our approach,
especially when $k$ is large, which increases the range of 
oscillations.
\smallskip

\subsection{Brief history}
Artin's primitive root conjecture states that given 
an integer $g$, possibly negative but not a perfect
square in $\Z$,  
the number $\p_N(g)$ of prime $p$ among 
$p_1=2,p_2=3,\ldots, p_N$ such that $g$ is primitive 
modulo $p$ approaches asymptotically $A(g)N$ as
$N\to \infty$ for
\begin{align}\label{amain}
&A(g)=A_h\, A_d,\ 
A_d=\Bigl(1-\,\mu(|d|)\prod_{p|d,\,p|h}\frac{1}{p-2}
\,\prod_{p|d,\,p\,\nmid\, h}\frac{1}{p^2-p-1}\Bigr),\\
&A_h\ =\ 
\prod_{p|h}(1-\frac{1}{p-1})\,
\prod_{p\,\nmid\, h}(1-\frac{1}{p(p-1)}), \hbox{ \ where}
\notag\\
d=&\hbox{\,Discriminant\,}(\Q[\sqrt{g}]),\ g=g_o^h
\hbox{ for } g_o\in \Z \hbox{ and maximal } h\in \N.
\notag  
\end{align}
Note that $\mu(|d|)= 0$ and, respectively,
$A(g)=A_h$ if and only if the discriminant $d$ is 
not from $1+4\Z$. If $h=1$ for such $g$, 
then $A(g)$ equals 
\begin{align}\label{aprod}
\hbox{Artin's constant }\ =\ 
A\ =\  \prod_{\hbox{\tiny prime\ } p}\,(1-\frac{1}{p(p-1)}).
\end{align}
According to \cite{St}, Artin's conjecture was finalized
around 1965. In 1967, it was deduced by
Hooley \cite{Ho} from the generalized Riemann hypothesis
for the fields $K_m=\Q[\ze_m,g^{1/m}]$ for squarefree
$m$. See \cite{Mo1} for a comprehensive introduction (including
some recent developments). See also \cite{Mu,Le}.  

Artin's heuristic approach to this conjecture was based on
the expectation that events ``{\em prime $p$ does not split completely 
in $K_q$ for prime $q$}" are independent (subject to later
qualitative and quantitative corrections). 
For instance, one can expect that $A$ equals 
$\lim_{N\to\infty} \p_N(g)/N$ if 
$g$ is {\em generic} as far as primitive roots
modulo prime $p$ are concerned. It leads to the following 
heuristic summation formula for Artin's constant:
\begin{align}\label{aadd}
A=\lim_{N\to \infty}\Si(N)\,,\ \Si(N)=\frac{1}{N}\sum_{i=1}^{N}
\frac{\phi(p_i-1)}{p_i-1}\,,
\end{align} 
which can be checked unconditionally, without any reference
to Artin's conjecture. See \cite{LL}, \cite{Mo1}. 
%
P. Moore extended it to
arbitrary $A(g)$  by switching to $p$
in this summation such that $g$ is a 
quadratic nonresidue modulo $p$ and $(p-1,h)=1$
(heuristically, it makes  
sense); see below.

As Lehmers wrote, the convergence in (\ref{aadd}) is 
``discouragingly slow" (they considered prime numbers $p<1500000$). 
It remains very slow when prime numbers
in much greater ranges are considered, generally,
no better than the (conjectural)
convergence of $\p_N(2)/N$ to $A=A(2)$; cf.
Table ``Artin's constant estimates" from
\cite{Si}\,($p<10^{14}$).   
\medskip

\subsection{Main Theorem}
A refine version of this heuristic approach is from
\cite{Mo1} (for any integer $g$):
\begin{align}\label{asumg}
&\lim_{N\to\infty}\frac{1}{N}\sum_{i=1}^N \varphi_g(p_i)= A(g)\ 
(\ \stackrel{\,\hbox{RH}}{\,=\kern -2pt =\,}\lim_{N\to\infty} 
\p_N(g)/N\,),
\where \\
&\varphi_g(p)\equal 2\,\frac{\phi(p-1)}{p-1} \for 
\Bigl(\frac{g}{p}\Bigr)=-1\and (p-1,h)=1,\notag \\  
&\varphi_g(p)\equal 0\, \hbox {\, otherwise\,.}\notag
\end{align} 

Heuristically, it is equally reasonable to expect that
\begin{align}\label{afrac}
A=\lim_{N\to \infty}\frac{\sum_{i=1}^{N}\phi(p_i-1)}
{\sum_{i=1}^N (p_i-1)}\,
\end{align} 
for sufficiently general $g$.
Switching here to the restricted summation from (\ref{asumg}),
we come to the following theorem.

\newtheorem{thno}{Theorem}

\begin{thno}\label{AG}
(i) For arbitrary integers $k\ge 0$ and $g$ (possibly,
negative),
\begin{align}\label{afracg}
&A(g)=\lim_{N\to \infty}\frac{\sum_{i=1}^{'\,N}p_i^k\phi(p_i-1)}
{\sum_{i=1}^{'\,N} p_i^k(p_i-1)},
\end{align} 
provided that $h=1$, where the summation $\Sigma'$ is over prime 
$p_i$ such that 
$\Bigl(\frac{g}{p_i}\Bigr)=-1$ and $(p_i-1,h)=1\,.$

(ii) If $h>1$, then the limit in the r.h.s.
of (\ref{afracg}) equals 
\begin{align}\label{afrach}
&\widetilde {A}(g) = \widetilde{A}_h\, \widetilde{A}_d\,,
\where \widetilde{A}_h\ =\  
\prod_{p\,\nmid\, h}(1-\frac{1}{p(p-1)})\,,\\
&\widetilde{A}_d\ =\ 1 \hbox{\, if \,} d\,|\,h \hbox{\ \ and\ \ }  
\widetilde{A}_d= A_d \hbox{\, otherwise\,.}\notag
\end{align}
\end{thno}
 
{\em Sketch of proof.} We follow \cite{LL,Mo1},
restricting ourselves with (\ref{afrac}); see also \cite{Pi}.
Coupling the generalized Landau formula from  \cite {SZ}
\begin{align}\label{landau} 
\sum_{p\le x}p^m=\frac{(1+o(1))\,x^{m+1}}{(1+m)\log x},\ m\ge 0,
\end{align}
with the classical estimate 
\begin{align}\label{pix}
\frac{\pi(x,d,1)}{\pi(x)}=\frac{1+O(1/\log x)}{\phi(d)}\,,
\end{align}
where $\pi(x,d,1)$ is the number of prime numbers $p\le x$
in $1+d\,\N$ ($\pi(x)=\pi(x,1,1)$), one arrives at:
\begin{align}\label{pix1}
\frac{\pi^{(m)}(x,d,1)}{\pi^{(m)}(x)}=
\frac{1+O(x^m/\log x)}{\phi(d)},
\end{align}
where $\pi^{(m)}(x,d,1)=\sum_{p\le x} p^m$  over 
prime $p\in 1+d\,\N$.
Then,
\begin{align*}
&\sum_{p\le x}\phi(p-1)=\sum_{p\le x} \sum_{d|p-1} 
(p-1)\frac{\mu(d)}{d}\\
=&\sum_{d|p-1}\frac{\mu(d)}{d}
\sum_{p\le x} (p-1)=\sum_{d\le x} 
\frac{\mu(d)}{d}(\pi^{(1)}(x,d,1)-\pi(x,d,1)).
\end{align*} 
Finally,
\begin{align*}
&\frac{\sum_{p\le x}\phi(p-1)}
{\sum_{p\le x} (p-1)}\,\thicksim\,\sum_{d\le x}\frac{\mu(d)}{d\phi(d)}
\,\stackrel{x\to \infty}{\,-\kern -5.5pt \to\,}\, A.
\end{align*}
We will omit the arguments (from the basic algebraic number theory)
that give the rational corrections for arbitrary $g$. \sq
 
\medskip
In (\ref{afrac},\ref{afracg}), only the leading powers of
$p_i$ matter; for instance, one can take here the ratio
$\sum_p\phi(\phi(p^{k+1}))/\sum_p\phi(p^{k+1})$, the
heuristic probability for $g$ being a primitive root
modulo $p^{k+1}$ over all prime $p$, which leads
to the same Artin constant.
\medskip

\subsection{Further examples} 
Let us apply the same procedure to
Stephens' constant $S(a,b)$ and the higher rank Artin
constants.  The former is defined for 
a given pair $a,b\in \Q^*$ such that 
$a^r b^s=1\Longrightarrow r=0,s=0$
for $r,s\in \Z$; it describes the density of prime $p$ such that 
$b=a^m\mod p$ for some $m\in \Z$. Modulo the generalized Riemann 
hypothesis, it equals
\begin{align}\label{steprod}
C_{ab}\,S \for S\equal\prod_{i=1}^\infty(1-\frac{p_i}{p_i^3-1}),
\end{align}
where the factors $C_{ab}$ are rational
\cite{S,MS}. These factors were calculated explicitly
in \cite{MS} under the condition that the group 
$\Q^*/\lan a,b,-1\ran$ is torsion free.
Assuming that $a,b$ are ``random", the
heuristic probability $P\!\!S(p)$ that $b=a^m \mod p$ for some $m$ 
can be readily calculated: 
\begin{align}\label{PS}
&P\!\!S(p)\ =\ \s(p)/(p-1)^2, \where 
\s(p)\equal \sum_{d|p-1}\,d\phi(d)\\
=&\prod_{j=1}^m
\frac{q_j^{2k_j+1}+1}{q_j+1}\ \hbox{\ for prime factorization\ }
\ p-1=\prod_{j=1}^m q_j^{k_j}.
\notag
\end{align}
We naturally omit the prime numbers that divide the numerators 
or denominators of $a,b$.
\smallskip

The rank $r$ Artin constant $A_r(g_1,\ldots, g_r)$
describes the heuristic density of prime\, $p$\,
such that a given set of nonzero integers 
$\{g_1,\ldots, g_r\}$ (or rationals) 
generates $\Z_p^*$\,. See \cite{CP}. 
Its ``generic" value (modulo 
the generalized Riemann hypothesis) is as follows:
\begin{align}\label{artr}
A_r\equal\prod_{i=1}^\infty\Bigl(\,1-\frac{1}{p_i^{\,r}\,(p_i-1)}\,
\Bigr).
\end{align}
The corresponding probability at $p$
for ``random" $g_1,\ldots,g_r$
equals
\begin{align}\label{PAr}
&P\!\!A_r(p)\ =\ \a_r(p)/(p-1)^r \for \a_r(p)\equal\prod_{j=1}^m
(q_j^{rk_j}- q_j^{r(k_j-1)})
\end{align}
in terms of the prime factorization $p-1=\prod_{j=1}^m q_j^{k_j}$.

The summation (unconditional) limiting formulas are as follows:
\begin{align}\label{sarlim}
&S\ =\ \lim_{N\to \infty}\frac{1}{N}\sum_{i=1}^{\,N} P\!\!S(p_i),\
\ A_r\ =\ \lim_{N\to \infty}\frac{1}{N}\sum_{i=1}^{\,N} 
P\!\!A_r(p_i).
\end{align} 
The proof of the following theorem is similar to that
of Theorem \ref{AG}. 

\begin{thno}\label{SA}
For an arbitrary integer $k\ge 0$,
\begin{align}\label{sform}
&S\ =\ \lim_{N\to \infty}\frac{\sum_{i=1}^{\,N} p_i^k\, \s(p_i)}
{\sum_{i=1}^{\,N} p_i^k\,(p_i-1)^2}\,,\\
&A_r\ =\ \lim_{N\to \infty}\frac{\sum_{i=1}^{\,N} 
p_i^k\, \a_r(p_i)}
{\sum_{i=1}^{\,N} p_i^k\,(p_i-1)^r}\,,\label{arform}
\end{align} 
where the summation is over consecutive
prime numbers $p_i$. \sq 
\end{thno}
\medskip

\subsection{Numerical aspects}
For the constants $C$ considered above,
we plot $\frac{\Si(N)}{C}-1$,
shown blue-thin, for the classical summation $\Si(N)$
and $\frac{R(N)}{C}-1$, which are red-thick, 
for our ratio approximations $R(N)$; 
for instance, $\Si(N)$ is from (\ref{aadd}) 
for Artin's $A$. The range is $N\le 1000$M. 



\vskip -.1in
\begin{figure}[htbp]
\begin{center}
\includegraphics[scale=1.2]{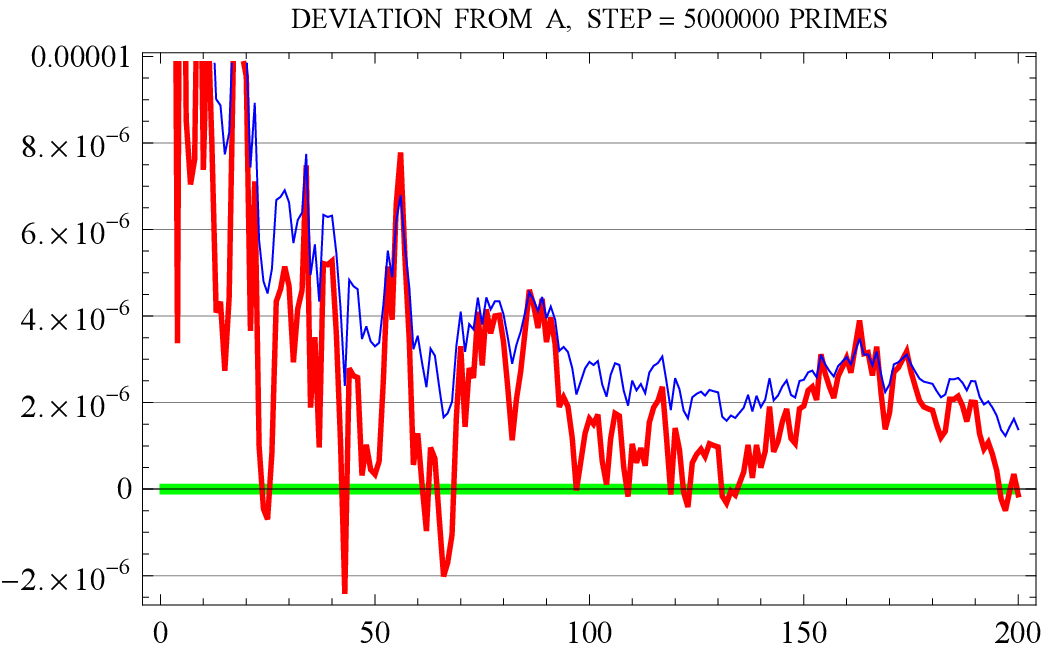}
\end{center}
\end{figure}

\vskip -0.35in
\begin{figure}[htbp]
\begin{center}
\includegraphics[scale=1.2]{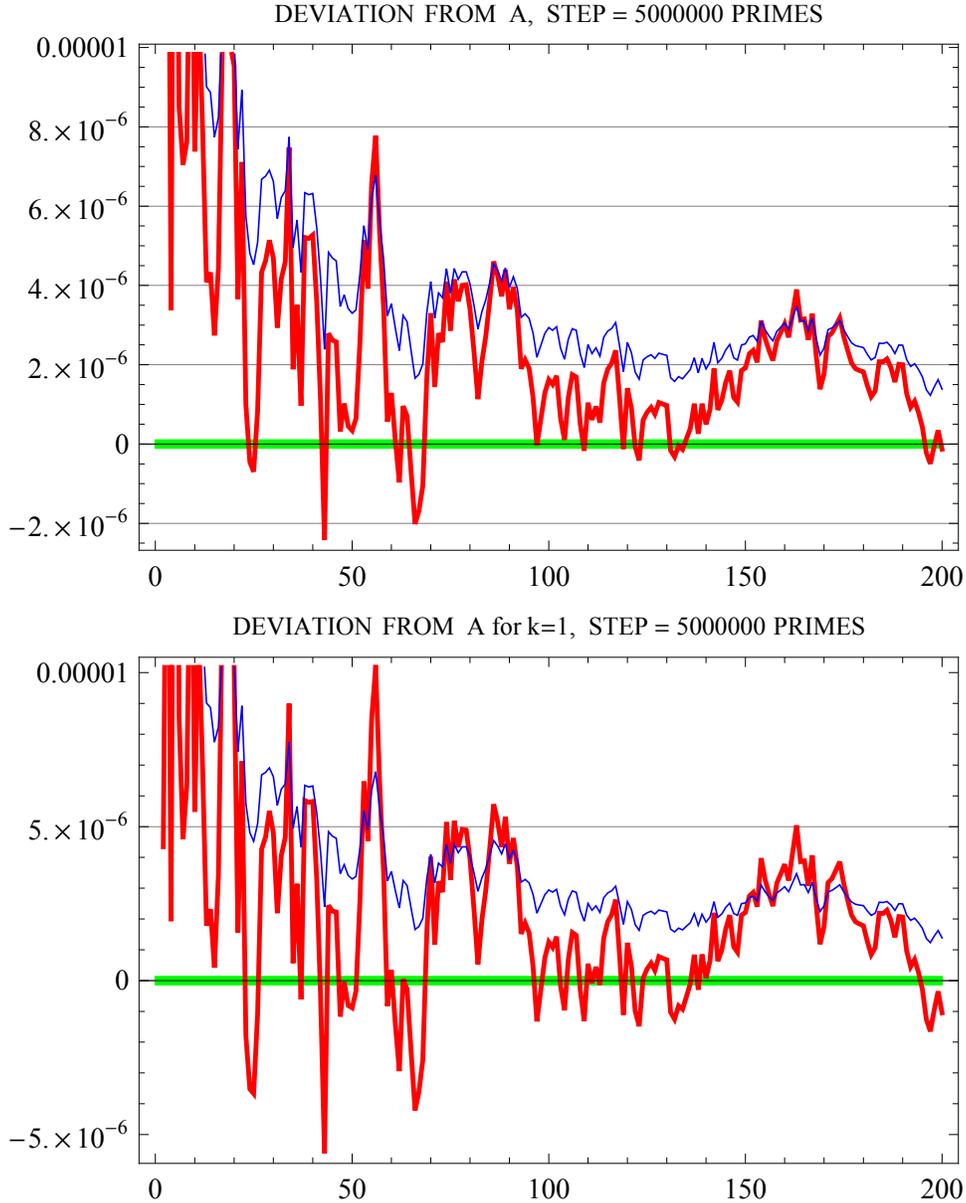}
\vskip -0.2in
\caption{Deviation from $A$ as $k=0,1$ for $1000$M primes. }
\label{art0-1}
\end{center}
\end{figure}
\smallskip

Figure \ref{art0-1} compares the stabilization 
of (\ref{aadd}), blue-thin, to the Artin constant
and the stabilization of 
(\ref{afrac}), red-thick, as $k=0,1$.
The function $\Sigma(N)-A$ remains
positive in the range $N\le 1000$M;
$R_k(N)$ oscillate
around Artin's constant $A\thickapprox 0.37395581361920228805$
(the zero level of this graph). The amplitude of oscillations
become larger for $R_{k=1}$ vs.  $R_{k=0}$,
however the graphs are very much similar.

We note that the best way to calculate $A$ 
and similar constants is based on the known 
product formulas in terms of 
$\ze(n)$ for integers $n>1$.

\begin{figure}[htbp]
\begin{center}
\includegraphics[scale=1.2]{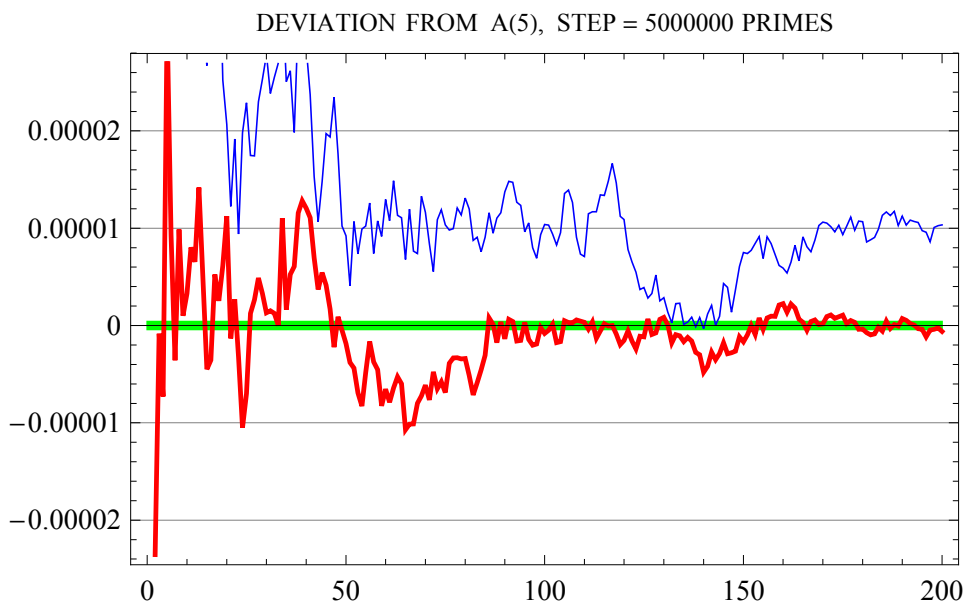}
\vskip -0.2in
\caption{Deviation from $A(5)$ for $1000$M primes. }
\label{art5}
\end{center}
\end{figure}

Figure \ref{art5}  shows the convergence
of (\ref{asumg}) and (\ref{afracg}) to $A(g=5)$ for $k=0$; 
notice that the blue-thin curve remains
beyond the red-thick one in this range. 
\medskip

 

\begin{figure}[htbp]
\begin{center}
\includegraphics[scale=1.2]{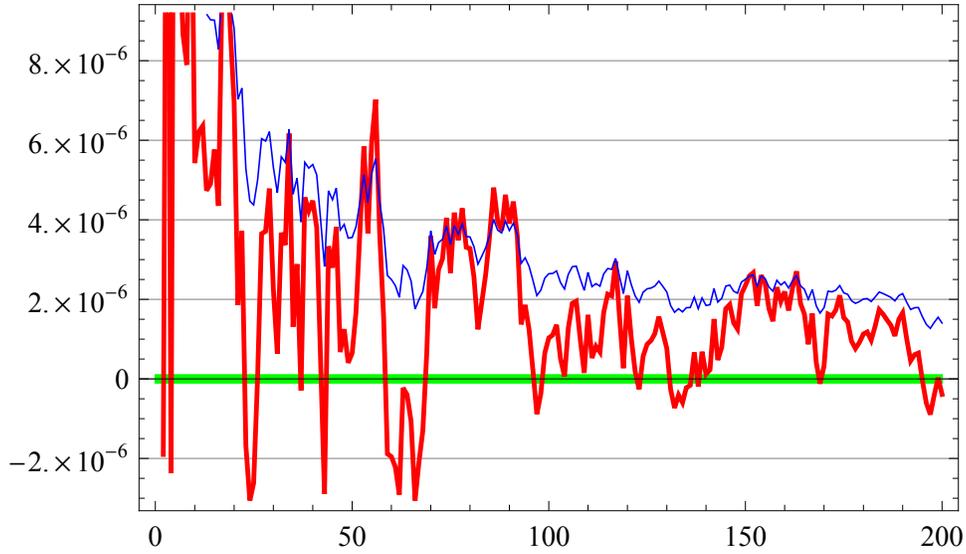}
\vskip -0.2in
\caption{Deviation from $S$  
for $1000$M primes. }
\label{steph}
\end{center}
\end{figure}

\begin{figure}[htbp]
\begin{center}
\includegraphics[scale=1.2]{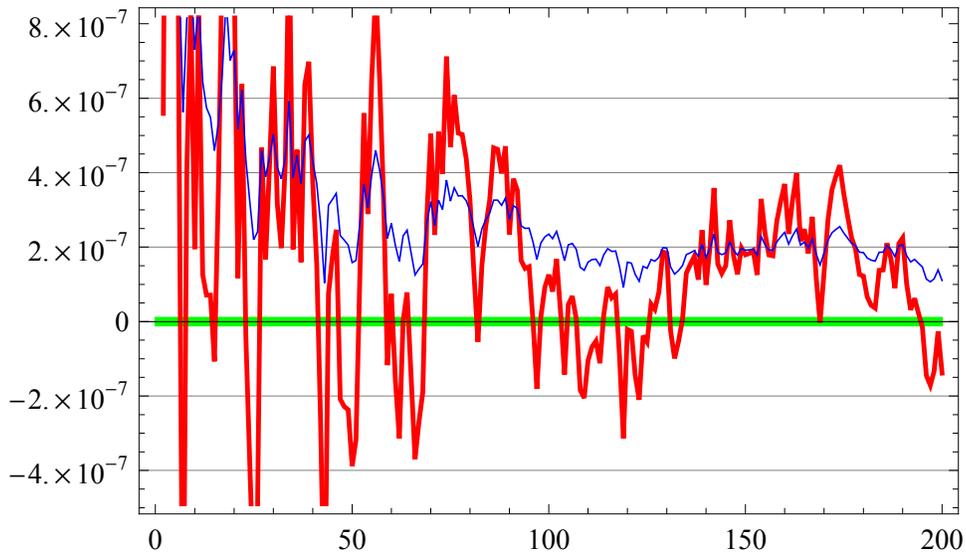}
\vskip -0.2in
\caption{Deviation from $A_3$\, (rank=$3$) 
for $1000$M primes. }
\label{3-art}
\end{center}
\end{figure}

The last two plots show the graphs for the Stephen constant
$S\thickapprox 0.57595996889294543964$, Figure \ref{steph}, 
and the rank 
$3$ Artin constant 
$A_3\thickapprox 0.85654044485354217443$, Figure \ref{3-art};
here $k=0$. 
The convergence rate and other features of these four graphs are 
similar to those for $A$ and $A(5)$. There is striking
(qualitative) similarity of these two figures, including the
oscillations, although the convergence 
rate in Figure \ref{3-art} (for $A_3$) is significantly
(almost $10$ times) greater than in Figure \ref{steph}. 
\medskip

Qualitatively, the behavior of $R(N)$ for large $N$ can
be evaluated following \cite{Pi}.
For instance, the functions $|R(N)/A-1|$ 
in Figure \ref{art0-1}  must be no greater than 
$C_{m}(\log N)^{-m}$ for {\em any} fixed $m>0$
and proper constant $C_{m}$ (depending on $k$) as $N>>0$.
Generally speaking, $C_m$ can be estimated in terms of 
(the order of) $N$, but we will not discuss it. Indeed, 
the graphs of $R(N)/A-1$, red-thick from Figure \ref{art0-1}, 
look like $O((\log N)^{-4})$ in the range $N<1000$M. 

The nature of oscillations of the functions $R(N)$ 
around the corresponding constants remains unclear.   
\medskip

{\em Acknowledgements}.
I am grateful to Zeev Rudnik for valuable comments.
I am very thankful to Pieter Moree for reading the
note, suggesting interesting questions toward
comparing our approach with the results from \cite{Mo2,Mo3}
and for the reference to \cite{Pi}.

\bibliographystyle{unsrt}

\end{document}

%% file: macro.tex
\newcommand{\Z}{{\mathbb Z}}
\newcommand{\Q}{{\mathbb Q}}
\newcommand{\N}{{\mathbb N}}

\def\HH{\mbox{${\mathcal H}$\kern-5.2pt${\mathcal H}$}}

\newtheorem{theorem }{Theorem}[section]
\newtheorem{maintheorem }[theorem]{Main Theorem}
\newtheorem{proposition }[theorem]{Proposition}
\newtheorem{definition }[theorem]{Definition}
\newtheorem{lemma }[theorem]{Lemma}
\newtheorem{corollary }[theorem]{Corollary}
\newtheorem{notation }[theorem]{Notation}
\newtheorem{remark }[theorem]{Remark}
\newtheorem{example }[theorem]{Example}

\newtheorem{ maintheorem }[theorem]{Main Theorem}
\newtheorem{ theorem}{Theorem}[section]
\newtheorem{ proposition}[theorem]{Proposition}
\newtheorem{ definition}[theorem]{Definition}
\newtheorem{ lemma}[theorem]{Lemma}
\newtheorem{ corollary}[theorem]{Corollary}
\newtheorem{ notation}[theorem]{Notation}
\newtheorem{ remark}[theorem]{Remark}
\newtheorem{ example}[theorem]{Example}

\def\for{\  \hbox{ for } \ }

\def\where{\  \hbox{ where } \ }
\def\and{\  \hbox{ and } \ }
\def\and{\  \hbox{ and } \ }

\def\equal{\stackrel{\,\mathbf{def}}{= \kern-3pt =}}

\def\Si{\Sigma}

\def\ze{\zeta}

\def\0{\mathbf{0}}

\def\çF{\mathcal{F}}

\def\p{\mathcal{P}}

\def\a{\mathcal{A}}

\def\s{\mathcal{S}}

\def\lan{\langle}

\def\ran{\rangle}

\newcommand{\sq}{\phantom{1}\hfill$\qed$}

\def\HH{\mathfrak{H}}

\def\HH{\hbox{${\mathcal H}$\kern-5.2pt${\mathcal H}$}}

\def\#{\sharp}
